\documentclass[12pt]{article}
\usepackage{amsmath,amsthm,amsfonts,amscd,amssymb,eucal,bbm}
\begin{document}
\def\RR{{\mathbb R}}
\def\CC{{\mathbb C}}
\def\NN{{\mathbb N}}
\def\ZZ{{\mathbb Z}}
\def\s1{{S^1}}
\def\diff{{{\rm Diff}^+(S^1)}}
\def\vect{{{\rm Vect}(S^1)}}
\def\sl2{{{\rm SL}(2,\RR)}}
\def\psl2{{{\rm PSL}(2,\RR)}}
\def\mob{{\rm Mob}}
\def\su2{{{\rm SU}(2)}}
\def\suN{{{\rm SU}(N)}}
\def\so3{{{\rm SO}(3)}}
\def\gk{{G_k}}
\def\nphi{{\mathfrak N}_\varphi}
\def\mphi{{\mathfrak M}_\varphi}
\def\aphi{{\mathfrak A}_\varphi}

\def\A{{\mathcal A}}
\def\B{{\mathcal B}}
\def\C{{\mathcal C}}
\def\D{{\mathcal D}}
\def\F{{\mathcal F}}
\def\H{{\mathcal H}}
\def\I{{\mathcal I}}
\def\K{{\mathcal K}}
\def\k{{\rm K}}
\def\M{{\mathcal M}}
\def\N{{\mathcal N}}
\def\O{{\mathcal O}}
\def\P{{\mathcal P}}
\def\R{{\mathcal R}}
\def\T{{\mathcal T}}
\def\U{{\mathcal U}}
\def\V{{\mathcal V}}
\def\W{{\mathcal W}}
\def\G{{\bf G}}

\newtheorem{theorem}{Theorem}[section]
\newtheorem{definition}[theorem]{Definition}
\newtheorem{corollary}[theorem]{Corollary}
\newtheorem{proposition}[theorem]{Proposition}
\newtheorem{lemma}[theorem]{Lemma}
\newtheorem{remark}[theorem]{Remark}
\title{{ INTERSECTING JONES PROJECTIONS}}

\author{SEBASTIANO CARPI
\\ Dipartimento di Scienze
\\ Universit\`a ``G. d'Annunzio" di Chieti-Pescara
\\ Viale Pindaro 87, I-65127 Pescara, Italy
\\ E-mail: carpi@sci.unich.it}
\date{} 
\maketitle 

\begin{abstract} Let $M$ be a von Neumann algebra on a Hilbert space 
$\H$ with a cyclic and separating unit vector $\Omega$ and let $\omega$ 
be the faithful normal state on $M$ given by 
$\omega(\cdot)=(\Omega,\cdot\Omega)$. 
Moreover, let  $\{N_i :i\in I\}$ be a family of von Neumann subalgebras 
of $M$ with faithful normal conditional expectations $E_i$ of $M$ onto 
$N_i$ satisfying  $\omega=\omega\circ E_i$ for all $i\in I$ and let 
$N=\bigcap_{i\in I} N_i$. We show that the projections $e_i$, $e$ 
of $\H$ onto the closed subspaces $\overline{N_i\Omega}$ and 
$\overline{N\Omega}$ respectively satisfy $e=\bigwedge_{i\in I}e_i$.
This proves a conjecture of V.F.R. Jones and F. Xu in \cite{JonesXu04}.
\end{abstract}

\section{Introduction} 
Let $M$ be a von Neumann algebra on a Hilbert space $\H$ and let 
$\Omega \in \H$ be a vector of norm $1$ which is cyclic and 
separating for $M$. Given a family $\{N_i: i\in I\}$ of von Neumann
subalgebras of $M$ it is often useful to consider the closed subspaces 
$\overline{N_i \Omega},\;i\in I$ and the corresponding projections 
$e_i\in {\bf B}(\H)$. If $N$ denotes the intersection 
$\bigcap_{i\in I}N_i$ and $e$ is the projection of $\H$ onto 
$\overline{N\Omega}$ one always have 
$e \leq \bigwedge_{i\in I} e_i$ namely 
$e\H \subset \bigcap_{i\in I} e_i\H$. However in general the equality 
does not hold and in fact it is not hard to give examples with 
$N=\CC1$ but $e_i =1$ for all $i\in I$ even when the set $I$ contains 
just two elements. 

In this paper we prove (see Corollary \ref{finalcorollary}) 
that if for every $i\in I$ there is a 
faithful normal conditional expectation of $M$ onto $N_i$ with 
$\omega \circ E_i =\omega$, where $\omega$ denotes the faithful 
normal state on $M$ given by $\omega (\cdot)=(\Omega, \cdot \Omega)$, 
then $e = \bigwedge_{i\in I}e_i$.

In a recent paper V.F.R. Jones and F. Xu gave a proof of 
this equality for a relevant class of examples associated with inclusions
of loop groups models of conformal nets on $\s1$ \cite[Lemma 4.14]{JonesXu04}. 
Moreover, they conjectured that this conclusion is not 
restricted to specific models but holds in general for inclusions of  
of completely rational conformal nets \cite[Remark 4.15]{JonesXu04}. 
In these examples $M$ is a local algebra of the larger conformal net
(a type ${\rm III}_1$ factor) and $\Omega$ is the vacuum vector
(if $M$ is a Type ${\rm II}_1$ factor and the vector state $\omega$ is the  
trace on M the equality holds by \cite{Skau}, cf. 
\cite[Remark 4.16]{JonesXu04}).

Our result proves the conjecture of Jones and Xu in \cite{JonesXu04} 
and actually shows that assumption (2) in \cite[Corollary 4.9]{JonesXu04}
is not needed. 

The proof we shall give is partially inspired by the one given in 
\cite[Lemma 4.14]{JonesXu04}. The main new idea is to replace the 
smeared vertex operators used in \cite{JonesXu04} by suitable 
closed operators provided by the Tomita-Takesaki modular theory. 
In fact we shall prove a more general result (Theorem 
\ref{maintheorem}) where a normal semifinite normal weight on the 
von Neumann algebra $M$ is given instead of the vector state $\omega$.

\section{Preliminaries and notations}
\label{definitionSect}
Let $M$ be a von Neumann algebra and let $\varphi$ be a normal 
semifinite faithful (n.s.f.) weight on $M$. Then the set 
\begin{equation}
\nphi =\{ x \in M : \varphi (x^*x) < \infty \}. 
\end{equation}
is $\sigma$-weakly dense left ideal of $M$. $\nphi$ with the inner 
product $(x,y)= \varphi(x^*y)$ can be completed to a complex Hilbert
space $\H_\varphi$. Accordingly $\nphi \subset M$ will be considered as 
a dense subspace of $\H_\varphi$ via the mapping 
$\nphi \ni x \mapsto x_\varphi \in \H_\varphi$ so that for 
$x,y \in \nphi$ we have $(x_\varphi, y_\varphi)= \varphi(x^*y)$. 
The GNS representation $\pi_\varphi$ of $M$ on $\H_\varphi$ is 
determined by $\pi_\varphi(x)y_\varphi = (xy)_\varphi$ for all 
$x\in M$ and $y\in \nphi$. Then $\pi_\varphi$ is normal and faithful i.e. 
a $\ast$-isomorphism of $M$ onto the von Neumann algebra 
$\pi_\varphi(M) \subset {\bf B}(\H_\varphi)$ (see \cite[\S 2]{Stratila}).    

The set $\nphi \cap \nphi^*$ is a $\sigma$-weakly dense 
self-adjoint subalgebra
$\aphi$ of $M$ which can also be considered as a dense subspace of 
$\H_\varphi$.  
The antilinear operator $S_\varphi^0$ on $\H_\varphi$ with domain $\aphi$, 
defined by
\begin{equation}
S_\varphi^0 x_\varphi = (x^*)_\varphi,\quad x\in \aphi, 
\end{equation}
is preclosed and we denote its closure by $S_\varphi$. The modular 
operator $\Delta_\varphi=S_\varphi^* S_\varphi$ and the modular 
conjugation $J_\varphi$ associated with $M$ and $\varphi$ are obtained 
from the polar decomposition $S_\varphi=J_\varphi\Delta_\varphi^{1/2}$ of 
$S_\varphi$. Moreover the following fundamental relations hold
\begin{equation}
J_\varphi \pi_\varphi (M)J_\varphi=\pi_\varphi(M)',\; 
\Delta_\varphi^{it}\pi_\varphi(M)\Delta_\varphi^{-it}=\pi_\varphi(M),\; 
t\in \RR 
\end{equation}
and the modular automorphism group $\{\sigma_t^\varphi \}_{t\in \RR}$ 
of $M$ associated with $\varphi$ is defined by 
\begin{equation}
\Delta_\varphi^{it}\pi_\varphi(x)\Delta_\varphi^{-it}
=\pi_\varphi(\sigma_t^\varphi (x)),\; x\in M, t\in\RR.
\end{equation}

For every $\eta \in \H_\varphi$ one defines a linear operator
$R_\eta^0$ on $\H_\varphi$ with domain $\aphi$ by
\begin{equation} 
R_\eta^0x_\varphi=\pi_\varphi(x)\eta,\; x\in \aphi.
\end{equation}
If $\eta$ is in the domain $D(S_\varphi^*)$ of $S_\varphi^*$ then 
$R_\eta^0$ is preclosed and its closure $R_\eta$ is affiliated 
with $\pi_\varphi(M)'$, see \cite[Chapter I, \S 2]{Stratila}. 
The subset $\aphi'\subset \H_\varphi$ defined by
\begin{equation}
\aphi'= \{\eta \in D(S_\varphi^*): R_\eta \in {\bf B}(\H_\varphi) \}
\end{equation}
is a dense subspace of $\H_\varphi$ and the set 
$\{R_\eta: \eta \in \aphi' \}$
is a $\sigma$-weakly dense self-adjoint subalgebra of 
$\pi_\varphi(M)'$. 

Similarly, for every $\xi \in \H_\varphi$ one defines a linear operator 
$L_\xi^0$ 
on $\H_\varphi$ with domain $\aphi'$ by 
\begin{equation}
\label{l0xi}
L_\xi^0\eta=R_\eta \xi,\; \eta \in \aphi'.
\end{equation}
If $\xi$ is in the domain $D(S_\varphi)$ of $S_\varphi$ then $L_\xi^0$ 
is preclosed and its closure $L_\xi$ is affiliated with 
$\pi_\varphi (M)$. 

The operators $L_\xi,\;\xi \in D(S_\varphi)$, which in general can be 
unbounded, will play a crucial role in the proof of our main result. 

We conclude this section with a proposition (cf. \cite{Skau}).
which we shall need later. 

\begin{proposition} 
\label{affiliatedprop}
Let $\{R_i: i\in I\}$ be a family of von Neumann algebras 
on a Hilbert space $\H$ and let $T$ be a closed linear operator on 
$\H$. If $T$ is affiliated with $R_i$ for every $i\in I$ then $T$ is also 
affiliated with $R=\bigcap_{i\in I}R_i$. 
\end{proposition} 
\begin{proof} Let $A$ be the unital self-adjoint subalgebra of 
${\bf B}(\H)$ generated by the union of the algebras $R_i', i\in I$. 
If $x\in A$ then $xT\subset Tx$. Now $A''=R'$ and thus $A$ is 
strong-operator dense in $R'$. For $x\in R'$ let $x_\lambda$ be a net
in $A$ converging to $x$ in the strong-operator topology. If $\xi$ is 
in the domain $D(T)$ of $T$ then $x_\lambda\xi \in D(T)$ and 
$Tx_\lambda \xi= x_\lambda T\xi$ for each $\lambda$. Hence 
$\lim x_\lambda \xi =x\xi$ and $\lim Tx_\lambda \xi =xT\xi$. 
Since $T$ is closed it follows that $x\xi \in D(T)$ and $Tx\xi=xT\xi$. 
Hence $xT\subset Tx$ for every $x\in R'$ namely $T$ is affiliated 
with $R$.
\end{proof} 

\section{Results}
Let $M$ be a von Neumann algebra and let $\varphi$ be a n.s.f. weight on 
$M$. If $N$ is a von Neumann subalgebra of $M$ we can define a closed 
subspace $\H_N$ of $\H_\varphi$ by
\begin{equation}
\H_N=\overline{\{x_\varphi: x\in \nphi \cap N\}}.
\end{equation}
For all $x\in N$ we have $x(\nphi \cap N) \subset \nphi \cap N$ and hence 
$\H_N$ is invariant for $\pi_\varphi(N)$.

If  $\sigma^\varphi_t(N)=N$ for all $t\in \RR$ and the restriction $\psi$ 
of $\varphi$ to $N$ is
semifinite we say that $N$ is a {\it modular covariant}
von Neumann subalgebra of $M$ relatively to $\varphi$
or simply a modular covariant subalgebra if the corresponding weight
on $M$ is unambiguously defined from the context.
Note that if $\varphi$ is a state, i.e. $\varphi(1)=1$, then $N\subset M$
is modular covariant iff $\sigma^\varphi_t(N)=N$ for all $t\in \RR$.

A von Neumann subalgebra $N\subset M$ is modular covariant if and only 
if there exists a faithful normal conditional expectation
$E$ of $M$ onto $N$ such that $\varphi = \varphi \circ E$, namely
$\varphi(x)= \varphi(E(x))$ for each positive element $x$ of $M$
\cite{Takesaki72} (see also \cite[\S 10]{Stratila}).
In this case the conditional expectation $E$ is completely determined 
by 
\begin{equation} 
\label{Ee}
\pi_\varphi(E(x))e=e\pi_\varphi(x)e,\; x\in M,
\end{equation}
where $e$ denotes the projection of $\H_\varphi$ onto 
$\H_N$ (the {\it Jones projection}), and the fact that $\H_N$ 
is separating for $\pi_\varphi(M)$,  being 
$N\cap\nphi$ $\sigma$-weakly dense in $N$. 

\begin{lemma}
\label{eprime}
 Let $\varphi$ be a n.s.f. weight on the von Neumann 
algebra $M$. If $N\subset M$ is a modular covariant von Neumann 
subalgebra and $e$ is the corresponding Jones projection then 
$$\pi_\varphi(N)=\pi_\varphi(M)\cap \{e\}'.$$
\end{lemma}
\begin{proof} From the fact that $\H_N$ is invariant for 
$\pi_\varphi(N)$ it follows that $e\in \pi_\varphi(N)'$. 
Assume now that $x\in M$ and that $\pi_\varphi(x)$ commutes 
with $e$. Then it follows from Eq. (\ref{Ee}) that 
$\pi_\varphi(E(x)-x)e=0$ and hence, being $\H_N$ separating 
for $\pi_\varphi(M)$ and $\pi_\varphi$ faithful, that 
$x=E(x)\in N$.
\end{proof}

\begin{proposition}
\label{lxiprop}
Let $M$, $\varphi$ and $N$ be as in the previous lemma 
and for $\xi\in D(S_\varphi)$ let $L_\xi$ be the closed operator 
affiliated with $\pi_\varphi(M)$ defined after Eq. (\ref{l0xi}). 
Then $L_\xi$ is 
affiliated with $\pi_\varphi(N)$ for all $\xi\in D(S_\varphi) \cap \H_N$.
\end{proposition} 
\begin{proof} Let $\xi\in D(S_\varphi) \cap \H_N$. Since we know that 
$L_\xi$ is a closed operator affiliated with  $\pi_\varphi(M)$ and 
by Lemma \ref{eprime} $\pi_\varphi(N)=\pi_\varphi(M)\cap \{e\}'$ it 
follows from Proposition \ref{affiliatedprop} that it is enough to 
show that $L_\xi$ is affiliated with $\{e\}'$ namely that 
$eL_\xi \subset L_\xi e$. 

From $eJ_\varphi = J_\varphi e$ (see pag. 131 of \cite{Stratila})  
and $\aphi'=J_\varphi \aphi$ \cite[2.12]{Stratila} 
it follows that $e\aphi'=J_\varphi e\aphi$. Now, for every 
$x\in \nphi$ we have $E(x)\in \nphi$ and $ex_\varphi =E(x)_\varphi$
, see \cite[10.3]{Stratila}. Since $E$ is a self-adjoint map it follows 
that $e\aphi \subset \aphi$ and hence that $e\aphi'\subset \aphi'$.  
Given $\eta\in \aphi'$ , $x\in \aphi$ we have 
\begin{eqnarray*}
eR_\eta e x_\varphi &=& eR_\eta E(x)_\varphi = e \pi_\varphi(E(x)) \eta 
 = \pi_\varphi(E(x))e \eta \\
&=& R_{e \eta}E(x)_\varphi = R_{e \eta}e x_\varphi. 
\end{eqnarray*}  
Since $R_\eta$ and $R_{e\eta}$ are bounded and $\aphi$ is dense in
$\H_\varphi$ it follows that $eR_{\eta}e=R_{e\eta} e$ 
for every $\eta\in \aphi'$. Hence, using the assumption that $\xi= e\xi$, 
for every $\eta \in \aphi'$ we find 
\begin{eqnarray*}
L_\xi e \eta &=& R_{e\eta}\xi = R_{e \eta} e \xi = e R_\eta e \xi \\
&=& eR_\eta \xi = e L_\xi \eta 
\end{eqnarray*} 
and the conclusion follows from the fact that $\aphi'$ is a core for 
$L_\xi$. 
\end{proof}

We are now ready prove the main result of this paper. 

\begin{theorem} 
\label{maintheorem}
Let $M$ be a von Neumann algebra with a n.s.f. weight 
$\varphi$ and let $\{N_i : i\in I \}$ be a family of modular 
covariant von Neumann subalgebras of $M$ with Jones projections 
$\{e_i: i\in I \}$. Assume that the restriction of $\varphi$ to 
$N= \bigcap_{i\in I}N_i$ is semifinite. Then $N$ is a modular covariant 
subalgebra of $M$ with Jones projection $e$ satisfying
$e=\bigwedge_{i \in I}e_i.$
\end{theorem}  
\begin{proof} We have to show that $\bigcap_{i\in I} \H_{N_i}=\H_N$. 
For all $i\in I$ we have $\nphi \cap N \subset \nphi \cap N_i$ and 
hence $\H_N \subset \bigcap_{i\in I}\H_{N_i}$. To prove the other 
inclusion let us consider the projection $f$ of $\H_\varphi$ onto 
$\bigcap_{i\in I} \H_{N_i}$. 
Since $e_i \Delta_\varphi \subset \Delta_\varphi e_i$ 
(see pag. 131 of \cite{Stratila}), $\Delta_\varphi$ is affiliated with 
$\{e_i\}'$ for each $i\in I$ and hence, by Proposition 
\ref{affiliatedprop} it is affiliated with 
$\bigcap_{i\in I}\{e_i\}' \subset \{f \}'$. 
It follows that 
$f \Delta_\varphi^{1/2} \subset \Delta_\varphi^{1/2} f$
and thus that $D(S_\varphi)\cap f\H_\varphi =fD(S_\varphi)$. 

Now let  $\xi \in D(S_\varphi)$. Then 
$f\xi \in \bigcap_{i\in I}\left( \H_{N_i}\cap D(S_\varphi) \right)$ and 
by Propositions \ref{affiliatedprop} and \ref{lxiprop} $L_{f\xi}$ is 
affiliated with $\bigcap_{i\in I} N_i = N$. It follows that 
$e L_{f\xi} \subset L_{f\xi} e$ and hence that 
$e R_\eta f\xi = R_{e\eta} f\xi$ for every $\eta\in \aphi'$. Thus, using 
the fact that  $\aphi'=J_\varphi \aphi$, we find
$R_{eJ_\varphi x_\varphi}f\xi\in e\H_\varphi$ for all $x\in \aphi$.
From the equalities 
$R_{J_\varphi y_\varphi}=J_\varphi\pi_\varphi(y)J_\varphi$, 
$y\in \aphi$ (see pag. 26 of \cite{Stratila}) and 
$J_\varphi e = e J_\varphi$ it follows that 
$J_\varphi\pi_\varphi(E(x))J_\varphi f\xi \in e\H_\varphi$ for every 
$x\in \aphi$, where $E$ is the faithful normal conditional expectation of 
$M$ onto $N$ satisfying $\varphi \circ E =\varphi$. Hence, being $\aphi$ 
$\sigma$-weakly dense in 
$M$ and $E$ normal, we find $f\xi \in e\H_\varphi =\H_N$. Since $\xi \in 
D(S_\varphi)$ 
was arbitrary we can conclude that $f\H_\varphi \subset \H_N$. 
\end{proof}

\begin{corollary}
\label{finalcorollary}
Let $M$ be a von Neumann algebra on a Hilbert space
$\H$ with a cyclic and separating unit vector $\Omega$ and let $\omega$
be the faithful normal state on $M$ given by 
$\omega(\cdot)=(\Omega,\cdot\Omega)$. Assume that for a given family 
$\{N_i :i\in I\}$ of von Neumann subalgebras of $M$ there exist 
faithful normal conditional expectations $E_i$ of $M$ onto $N_i$ 
satisfying $\omega=\omega\circ E_i$ for all $i\in I$ and let
$N=\bigcap_{i\in I} N_i$. Then the Jones projections $e_i$, $e$
of  $\H$ onto the closed subspaces $\overline{N_i\Omega}$ and
$\overline{N\Omega}$ respectively satisfy $e=\bigwedge_{i\in I}e_i$.
\end{corollary}

\bigskip 

\noindent{\bf Acknowledgements} 
The author would like to thank R. Longo for discussions.
He also thanks Y. Kawahigashi for the invitation and the 
hospitality at the Department of Mathematical Sciences of the 
University of Tokyo in October 2004 where part of this work was 
done.

\end{document}